\documentclass[a4paper,11pt,oneside]{article}

\usepackage{amsmath,amsthm,amsfonts,amssymb,amscd}
\usepackage{subfigure}
\usepackage{graphicx,xcolor}
\usepackage{algorithm}
\usepackage[ngerman,english]{babel}

\usepackage{units}  
\usepackage{siunitx}

\usepackage{authblk}

\usepackage[top=2cm, bottom=2cm, left=2cm, right=2cm]{geometry}

\theoremstyle{plain} 
\newtheorem{theorem}{Theorem}[section]

\newtheorem{Remark}[theorem]{Remark}
\newtheorem{Definition}[theorem]{Definition}

\newtheorem{form}[theorem]{Formulation}

\theoremstyle{definition} %

\newtheorem{algo}{Algorithm}

\theoremstyle{remark} %

\usepackage{xcolor}

\newcommand{\bs}[1]{{#1}} 

\newcommand{\hOmega}{\widehat\Omega}

\newcommand{\hGamma}{\widehat\Gamma}
\newcommand{\hI}{\widehat I}

\newcommand{\R}{\mathbb R}
\newcommand{\hX}{\bs{\widehat{X}}}

\newcommand{\hV}{\bs{\widehat{V}}}

\newcommand{\hL}{{\widehat{L}}}

\newcommand{\IT}{I}

\newcommand{\hcalA}{\widehat{\cal A}}

\newcommand{\hF}{\bs{\widehat{F}}}
\newcommand{\hJ}{\widehat{J}}
\newcommand{\hE}{\bs{\widehat{E}}}

\newcommand{\hdiv}{\widehat{\operatorname{div}}}

\newcommand{\hnabla}{\bs{\widehat\nabla}}
\newcommand{\hSigma}{\bs{\widehat\Sigma}}

\newcommand{\hsigma}{\bs{\widehat\sigma}}

\newcommand{\hrho}{\bs{\widehat\rho}}

\newcommand{\hx}{\bs{\widehat x}}

\newcommand{\hv}{\bs{\widehat v}}
\newcommand{\hu}{\bs{\widehat u}}
\newcommand{\hp}{\bs{\widehat p}}

\newcommand{\hw}{\bs{\widehat w}}

\newcommand{\hg}{\bs{\widehat g}}

\newcommand{\hA}{\bs{\widehat A}}

\newcommand{\hU}{\bs{\widehat U}}

\newcommand{\hZ}{\bs{\widehat Z}}

\newcommand{\hPsi}{\bs{\widehat\Psi}}

\newcommand{\hpsi}{\bs{\widehat\psi}}

\newcommand{\hn}{\bs{\widehat n}}

\newcommand{\pt}{\partial_t}

\begin{document}

\title{
Optimization with nonstationary, nonlinear monolithic fluid-structure interaction
}
\author[1]{Thomas Wick}
\author[2]{Winnifried Wollner}

\affil[1]{Leibniz Universit\"at Hannover, Institut f\"ur Angewandte
  Mathematik, AG Wissenschaftliches Rechnen, Welfengarten 1, 30167 Hannover, Germany}
\affil[2]{Technische Universit{\"a}t Darmstadt, Fachbereich Mathematik,
Dolivostrasse 15, 64293~Darmstadt, Germany}

\date{}

\maketitle
	
\begin{abstract}
Within this work, we consider optimization settings for 
nonlinear, nonstationary fluid-structure interaction. 
The problem is formulated in a monolithic fashion using 
the arbitrary Lagrangian-Eulerian framework to set-up
the fluid-structure forward problem. 
In the optimization approach,
either optimal control or parameter estimation
problems are treated. In the latter, 
the stiffness of the solid is estimated from given reference values. In the numerical 
solution, the optimization problem is solved
with a gradient-based solution algorithm. The nonlinear subproblems of the FSI forward problem  are solved with a Newton method 
including line search. Specifically, we will formally provide the
backward-in-time running adjoint state used for gradient computations. Our algorithmic 
developments are demonstrated with some numerical examples as for instance 
extensions of the well-known fluid-structure benchmark settings
and a flapping membrane test in a channel flow with elastic walls.
\end{abstract}

\section{Introduction}
\label{sec_intro}
This paper is devoted to the study of 
optimal control and parameter estimation problems 
of nonstationary, nonlinear fluid-structure interaction (FSI).
For general overviews on the FSI forward problem, we refer to the books 
\cite{BuSc06,FoQuaVe09,GaRa10,BuSc10,BaTaTe13,BoGaNe14,Ri17_fsi,FrHoRiWiYa17}.
Fluid-structure interaction 
is still one of the most 
challenging problem settings within multiphysics 
applications.
The main reason being that the dynamics 
of both subproblems are exchanged on the interface and 
accurate discretizations are therefore necessary. 
Secondly, numerical algorithms are sensitive in terms 
of stability to the physical parameters; known as added-mass 
effect~\cite{CaGeNo05,FoeWaRa07,AstChoFe09,Brumm09}.
As coupling strategy we choose, in this work, the 
well-known arbitrary Lagrangian-Eulerian 
(ALE) technique~\cite{DoFaGi77,HuLiZi81,FoNo99}.

Employing FSI as forward problem within an optimization framework 
contains the previously mentioned difficulties and yields significant
further challenges when dealing with nonstationary problem settings.
Historically, this subject falls into the category of PDE-constrained 
optimization~\cite{Lions:1971}.
Studies concentrating on theoretical and computational aspects for stationary FSI optimization are
\cite{RiWi13_fsi_opt,Wi11_phd,WiWo19,Chirco_2017}.
Here, we notice that the required adjoints are the same 
as used for adjoint-based error estimation; see for instance
\cite{ZeeBrummelenAkkermanBorst2011,Ri12_dwr}.
Nonlinear (stationary) FSI investigating various partitioned coupling 
techniques was recently subject in~\cite{Sing19}.
The by far more challenging situation of nonstationary settings 
is listed in the following.  
A nonstationary situation assuming a rigid solid was theoretically 
studied in~\cite{MouZol06}. 
Further theoretical results for 
a boundary control FSI problem were established in~\cite{BucciLasiecka:2010}.
Parameter estimation to detect the stiffness of an arterial 
wall with a well-posedness analysis and numerical simulations 
was addressed in~\cite{PeVeVe11}. Again in blood flow simulations, 
a data assimilation problem was formulated in~\cite{GUERRA201457},
in which however, the arterial walls were not considered.
A full FSI problem for data assimilation using a Kalman filter 
was subject in~\cite{BeMoiGer12}.
In~\cite{Kuberry2013594}, the authors used optimization techniques 
to formulate the FSI coupling conditions.
Adjoints for 1D FSI were derived in~\cite{Degroote2013,Maritin2005}.
Reduced basis methods for FSI-based optimization were developed
in~\cite{Lassila2013}.
Optimal control of nonstationary FSI applied to benchmark settings 
was investigated in \cite{BAZILEVS20131989}.
A linearized FSI optimization problem 
was addressed in~\cite{FaiMeiVex16} and detailed results 
for full-time-dependent FSI optimal control were summarized in~\cite{Fai17}.
In this respect, we also mention~\cite{FaiWi18} in which the 
adjoints required for optimization were employed for 
dual-weighted residual error estimation for time adaptivity.
Most recently, a uncertainty 
quantification framework for fluid-structure interaction with 
applications in aortic biomechanics was 
developed in~\cite{Kratzke18}.

The significance of the current work is on the development 
of a robust fully monolithic formulation for gradient-based optimization
for nonstationary, nonlinear FSI problems. Here, the coupled 
problem is prescribed in the reference configuration 
with the help of the ALE approach in a variational-monolithic way.
As previously 
summarized in our literature review, only very few results exist to date
for such a framework. Indeed the challenges consist of both the nonlinearities
and the nonstationary nature of the problem. FSI in the forward 
state is itself a highly nonlinear problem. Moreover, interesting nonstationary 
configurations require several thousands of time steps. For instance
the FSI 3 benchmark~\cite{HrTu06b,BuSc10} requires about $6\,000$ to $10\, 000$ time 
steps for a fully developed oscillatory solution. These are costly 
computations, even for a moderate number of spatial degrees of freedom.
Numerically, an inf-sup stable spatial discretization is applied 
to the FSI forward problem. Time discretization is based 
on a one-step-theta formulation. The discretized subproblems 
are solved with Newton solver including line search.
In order to apply gradient-based techniques, the adjoint state 
is running backwards in time and must access the primal solution
at the time points when treating nonlinear problems. Such 
derivations and implementations are very tedious. In this work,
we carefully derive and implement them in order to test their 
performance. 
These are tested with the help of 
the modification of well-known FSI benchmark settings~\cite{HrTu06b,BuSc10}
and a flapping membranes example that was originally proposed in~\cite{GiCaBoHa10}
and later modified in~\cite{Wi12_fsi_eale_heart}.

The outline of this paper is as follows: In Section~\ref{sec_eq}, 
the equations for fluid flow and solids are summarized. Moreover, the 
FSI setting is formulated in a monolithic fashion using the 
arbitrary Lagrangian-Eulerian framework. Section~\ref{sec_dis}
contains temporal and spatial discretizations. The main 
results are presented in Section~\ref{sec_gradient} in which 
the gradient computation, including details on the adjoint, are presented.
In Section~\ref{sec_algo}
the solution algorithms for the FSI optimization framework are presented.
Our algorithmic techniques are substantiated with several 
numerical tests in Section~\ref{sec_tests}. We summarize our main 
findings in Section~\ref{sec_conclusions}. 




\section{Modeling the FSI forward problem}
\label{sec_eq}
\subsection{Notation}
We denote by $\Omega :=\Omega (t) \subset \mathbb{R}^d$, $d=2$, 
the domain of the 
FSI problem. 
The domain consists of two time-dependent subdomains 
$\Omega_f (t)$ and $\Omega_s (t)$. 
The FSI-interface between $\Omega_f (t)$ and $\Omega_s (t)$  is denoted by
$\Gamma_i (t) = \overline{\partial\Omega_f} (t) \cap 
\overline{\partial\Omega_s} (t)$.  
The initial (or later reference) 
domains are denoted by $\hOmega, \hOmega_f$ and $\hOmega_s$, respectively, with the
interface $\hGamma_i = \overline{\partial\hOmega_f} \cap \overline{\partial\hOmega_s}$.
Furthermore, we denote the outer boundary by 
$\partial\hOmega = \hGamma = \hGamma_{\text{in}} \cup \hGamma_D \cup \hGamma_{\text{out}}$
where
$\hGamma_D $ and $\hGamma_{\text{in}}$ are Dirichlet boundaries (for the velocities
and displacements) and 
$\hGamma_{\text{out}}$ denotes a fluid outflow Neumann
boundary, respectively. The displacements are set to zero on $\hGamma_{\text{out}}$.

As frequently used in the literature, 
we denote the $L^2$ scalar product in $\Omega$ 
with $(a,b) := (a,b)_{\Omega}:=\int_{\Omega}a\cdot b\,\mathrm{d}x$
for vectors $a,b$. For (2nd order) tensor-valued functions $A,B$, it yields
  $(A,B) := (A,B)_{\Omega}:=\int_{\Omega}A:B\,\mathrm{d}x$,
where $A:B = \sum_{ij=1}^d A_{ij}B_{ij}$ and $A_{ij}$ 
and $B_{ij}$ are the entries of $A$ and $B$.

\subsection{Spaces}
For the function spaces in the (fixed) reference domains $\hOmega, \hOmega_f,
\hOmega_s$, we define spaces 
for spatial discretization only. Rather 
than employing Bochner-spaces~\cite{DauLio2000,Wlo82} for space-time functions, the 
time $t$ is later explicitly accounted for, e.g.,~\cite{Evans:2000} (Section 7.1).
Here, let $\IT:=[0,T]$ be the time interval and $T$ the end time value.
First we define
\[
\hV:= H^1 (\hOmega)^d.
\]
Next, in the fluid domain, we define further:
\begin{align*}
\hL_f   &:= L^2 (\hOmega_f),\\
\hL_f^0 &:= L^2 (\hOmega_f)/\mathbb{R},\\
\hV_f^0 &:= \{ \hv_f\in H^1 (\hOmega_f)^d : \, 
\hv_f = 0 \text{ on } \hGamma_{\text{in}}\cup \hGamma_{D} \},\\
\hV_{f,\hu}^0 &:= \{ \hu_f\in H^1 (\hOmega_f)^d : \, 
\hu_f = \hu_s \text{ on } \hGamma_i, \quad \hu_f = 0 \text{ on } \hGamma_{\text{in}}
\cup \hGamma_{D}\cup\hGamma_{\text{out}} \}, \\
\hV_{f,\hu,\hGamma_i}^0 &:= \{ \hpsi_f\in H^1 (\hOmega_f)^d : 
\, \hpsi_f = 0 \text{ on } 
\hGamma_i \cup \hGamma_{\text{in}} \cup \hGamma_{D}\cup\hGamma_{\text{out}}\}.
\end{align*}
In the solid domain, we use
\begin{align*}
\hL_s   &:= L^2 (\hOmega_s)^d,\\
\hV_s^0 &:= \{ \hu_s\in H^1 (\hOmega_s)^d : \, 
\hu_s = 0 \text{ on } \hGamma_{D} \}.
\end{align*}

For the FSI problem using variational-monolithic coupling~\cite{HrTu06a,Du07,DuRaRi09} the velocity 
spaces are extended from $\hOmega_f$ 
and $\hOmega_s$ to the entire domain $\hOmega$ such that we can work with global $H^1$ functions.
Thus, we define:
\begin{equation}
\label{eq_global_H1_velocity_space}
\hV^0 := \{ \hv\in H^1 (\hOmega)^d : \, 
\hv = 0 \text{ on } \hGamma_{\text{in}}\cup \hGamma_{D} \}.
\end{equation}
By this choice, the kinematic and dynamic coupling conditions 
are automatically satisfied in a variational sense.

Finally, we notice that the spaces on the current domains $\Omega,\Omega_f, \Omega_s$ 
are defined correspondingly without `hat' notation.

\subsection{The ALE concept, transformed fluid flow, and solids in Lagrangian coordinates}
\label{section_ALE}
In this section, we recapitulate the ingredients to 
formulate a coupled problem (i.e., fluid-structure interaction)
 with the help of the ALE approach. The ALE mapping
 $\hcalA:\hOmega_f\to\Omega_f$ is defined first. 

\subsubsection{The ALE transformation and ALE time-derivative}
First, we define the ALE transformation:
\begin{Definition}
\label{def_ALE_mapping}
The ALE mapping is defined in terms of the vector-valued (artificial) fluid mesh displacement
$\hu_f:\hOmega_f\to\mathbb{R}^d$ such that 
\begin{equation}
\hcalA(\hx, t):\hOmega_f\times \IT \to\Omega_f, \quad\text{with }
\hcalA(\hx,t) = \hx+\hu_f(\hx,t),
\end{equation}
which is specified through
the deformation
gradient and its determinant
\begin{equation}
\hF:=\hnabla\hcalA = \hI +\hnabla\hu_f,\quad
\hJ:=\operatorname{det}(\hF).
\end{equation}
Furthermore, function values in Eulerian and Lagrangian coordinates
are identified by
\begin{equation}
u_f(x) =: \hu_f(\hx) , \quad \text{with } x = \hcalA(\hx,t).
\end{equation}
Here, $\hI$ denotes the identity matrix.
The mesh velocity is defined by $\hw := \partial_t\hcalA$.
The key quantity to measure the fluid mesh regularity is $\hJ$.
The artificial fluid displacement $\hu_f$ (the mesh motion) 
is obtained in this work by solving 
a biharmonic equation \cite{He01,Wi11,Du07,Wi11_phd}.

Finally, the transformation between different coordinate systems requires
transformation of derivatives. For a vector-valued function $u\in\Omega$ and $\hu\in\hOmega$
it holds, e.g.,~\cite{Ho00}:
\begin{equation*}
\nabla u = \hnabla\hu \hF^{-1}.
\end{equation*}
\end{Definition}
Finally, the ALE time-derivative is the total derivative 
of an Eulerian field and is important when working on moving domains:
\begin{equation}
\label{notations_ALE_time_derivative}
\partial_t |_{\hcalA} v_f(x,t) = \hw\cdot \nabla v_f + \partial_t v_f(x,t).
\end{equation}

\subsection{Equations for fluids and solids}
In this section, we briefly state the basic underlying equations 
first separately. 
In the following, we first present fluid flow and then the solid part.

{
\subsubsection{Strong forms}
\label{sec_strong_forms_FSI}
The isothermal, incompressible Navier-Stokes equations in an ALE
setting read: Given $v_{\text{in}}$, $h_f$, and $v_0$;
find $v_f:\Omega_f(t)\times I\to\mathbb{R}^d$ and 
$p_f:\Omega_f(t)\times I\to\mathbb{R}$ such that
\begin{align*}
&\rho_f \partial_t |_{\hcalA} v_f + \rho_f (v_f-\hw)\cdot\nabla v_f
-\nabla\cdot\sigma_f(v_f,p_f) = 0, \qquad \nabla\cdot v_f = 0
\quad\text{in }\Omega_f(t),\\
&v_f^D = v_{\text{in}}\text{ on } \Gamma_{\text{in}},\qquad
v_f = 0 \text{ on } \Gamma_D, \qquad
-p_fn_f + \rho_f\nu_f \nabla v_f\cdot n_f = 0 \text{ on } \Gamma_{\text{out}}, \qquad
v_f = h_f \text{ on } \Gamma_i,\\
&v_f(0) = v_0 \text{ in } \Omega_f(0),
\end{align*}
where the (symmetric) Cauchy stress is given by 
\[
\sigma_f(v_f,p_f) := -pI + \rho_f \nu_f (\nabla v + \nabla v^T),
\]
with the density $\rho_f$ and the kinematic viscosity $\nu_f$.
Later in the FSI problem, the function $h_f$ will 
be given by the solid velocity $v_s$. The normal vector 
is denoted by $n_f$.

The equations for geometrically non-linear elastodynamics are given 
as follows: Given $\widehat h_s$, $\widehat u_0$, and $\widehat v_0$;
find $\hu_s:\hOmega_s\times I\to \mathbb{R}^d$ such that
\begin{align*}
&\hrho_s \partial^2_t \hu_s - \hnabla\cdot(\hF\hSigma) = 0 \quad\text{in } \hOmega_s,\\
& \hu_s = 0 \text{ on }\hGamma_D, \qquad
\hF\hSigma \cdot \hn_s = \widehat h_s \text{ on }\hGamma_i,\\
&\hu_s(0) = \hu_0 \text{ in } \hOmega_s\times\{0\}, \qquad
\hv_s(0) = \hv_0 \text{ in } \hOmega_s\times\{0\}.
\end{align*}
The constitutive law is given by the tensor:
\begin{equation}\label{eq:Sigma}
\begin{aligned}
\hSigma &= \hSigma_s(\hu_s) = 2\mu \hE + \lambda tr(\hE) I, 
\quad\text{with } \hE = \frac{1}{2} (\hF^T \hF - I).
\end{aligned}
\end{equation}
Here, $\mu$ and $\lambda$ are the Lam\'e coefficients for the solid.
The solid density is denoted by $\hrho_s$ and the solid 
deformation gradient is $\hF = \hI + \hnabla\hu_s$.
Later in FSI, the vector-valued function $\widehat h_s$ will be given by the 
normal stress from the fluid problem.
Furthermore, $\hn_s$ denotes the normal vector.

}

\subsubsection{Variational forms}
The previous Navier-Stokes equations in a variational 
ALE framework described in 
a reference domain $\hOmega_f$ are given by:
\begin{form}[ALE Navier-Stokes in $\hOmega_f$]
\label{problem_NSE_ALE_fx}
Let $\hv_f^D$ a suitable extension of Dirichlet inflow data. 
Find vector-valued velocities and a scalar-valued pressure 
$\{\hv_f , \hp_f\} \in \{\hv_f^D + \hV_f^0\} \times \hL_f^0$ such
that the initial data $\hv_f (0) = \hv_f^0$ are satisfied, 
and for almost all times
$t\in \IT$ holds:
\begin{align*}
\hrho_f (\hJ\partial_t\hv_f , \hpsi_f^v)_{\hOmega_f} 
+ \hrho_f (\hJ\hF^{-1} (\hv_f - \hw)\cdot\hnabla\hv_f , 
\hpsi_f^v)_{\hOmega_f} 
+ (\hJ\hsigma_f \hF^{-T}, \hnabla\hpsi_f^v)_{\hOmega_f} \\
-\langle \hJ\hg_f \hF^{-T} \hn_f , \hpsi_f^v \rangle_{\hGamma_{\text{out}}} 
- \langle \hJ\hsigma_f \hF^{-T} \hn_f , \hpsi_f^v \rangle_{\hGamma_i} 
&= 0 \quad\forall\,\hpsi_f^v \in \hV_f^0,\\
(\hdiv\, (\hJ\hF^{-1} \hv_f) , \hpsi_f^p )_{\hOmega_f} 
&= 0\quad\forall\,\hpsi_f^p \in \hL^0_f.
\end{align*}
Here,
$\hg_f := - \hrho_f \nu_f \hF^{-T}\hnabla\hv_f^T$ denotes 
a correction term on the outflow boundary  and $\hn_f$ is the outer normal vector.
The transformed Cauchy stress tensor reads:
\begin{equation}
\label{transformed_Cauchy_stress_fluid}
\hsigma_f = 
-\hp_f\hI + 2\hrho_f\nu_f (\hnabla \hv_f\hF^{-1} + \hF^{-T} \hnabla \hv_f^T).
\end{equation}
\end{form}

The variational formulation for elastodynamics can be formulated 
as a first-order-in-time system:
\begin{form}[First order system in time weak formulation of elasticity
    including strong damping]
\label{problem_mixed_form_elasticity_with_damping}
Find $\hu_s\in \hV_s^0$ and $\hv_s\in \hL_s$ 
with the initial data $\hu_s(0) = \hu_0$ and $\hv_s(0) = \hv_0$
such that for almost all times $t\in\IT$:
\begin{align*}
\hrho_s(\pt \hv_s,\hpsi_s^v)_{\hOmega_s} 
+ (\hF\hSigma,\hnabla\hpsi_s^v)_{\hOmega_s} -\langle \hF\hSigma\hn_s , \hpsi_s^v \rangle_{\hGamma_{i}} 
&= 0
\quad\forall\, \hpsi_s^v\in \hV_s^0,\\
\hrho_s(\partial_t \hu_s - \hv_s, \hpsi_s^u)_{\hOmega_s} &= 0  \quad\forall \,\hpsi_s^u\in \hL_s.
\end{align*}
\end{form}

\subsection{Variational-monolithic ALE fluid-structure interaction}
\label{section_FSI}
\subsubsection{FSI interface coupling conditions}
\label{sec_interface_FSI}
The coupling of 
a fluid with a solid
must satisfy two physical conditions; namely  
continuity of velocities and continuity of normal stresses.
A third condition of geometric nature is necessary when working 
with the ALE framework: 
continuity of displacements, which couples the physical solid $\hu_s$ and
the fluid mesh motion $\hu_f$.
Mathematically, the first and third conditions can
be classified as (non-homogeneous) Dirichlet conditions and the second 
condition is a (non-homogeneous) Neumann condition.

In variational-monolithic coupling these Dirichlet conditions are built into
the corresponding function space by employing a 
globalized Sobolev space $\hV^0$ (see~\eqref{eq_global_H1_velocity_space}). 
Neumann type conditions are weakly 
incorporated through interface integrals, which 
actually cancel out in the later models because of their 
weak continuity thanks to working with the space $\hV^0$.

For the fluid problem, continuity of velocities is required (i.e., 
a kinematic coupling condition):
\begin{equation}
\label{coupling_velocities_strong_form}
\hv_f = \hv_s \quad\text{on } \hGamma_i .
\end{equation}
To complete the solid problem, we must enforce 
the balance of the 
normal stresses on the interface (i.e., a dynamic coupling condition):
\begin{equation}
\label{coupling_normal_stresses_strong_form}
\hJ\hsigma_f \hF^{-T}\hn_{f} + \hF\hSigma\hn_{s} = 0 \quad\text{on } \hGamma_i .
\end{equation}
For the geometric problem we have
\begin{equation}
\label{equations_definition_EXT}
\hu_f = \hu_s \quad\text{on } \hGamma_i,
\end{equation}
from which we obtain immediately $\partial_t \hu_s = \hv_s = \hv_f$ on $\hGamma_i$ 
by temporal differentiation.

\subsubsection{The FSI model using biharmonic mesh motion}
Combining the previous equations for fluids and solids and 
applying biharmonic mesh motion for realizing the 
ALE mapping, we obtain the following FSI model~\cite{Du07,Wi11_phd,Wi11}:

\begin{form}[Variational-monolithic ALE FSI in $\hOmega$]
\label{eq:fsi:ale:harmonic} 
 Let the constitutive laws from before be given and $\widehat \alpha >0$ be a small parameter.
Find a global vector-valued velocity, vector-valued displacements, additional 
displacements (due to the splitting of the biharmonic mesh motion model into 
two second-order equations)  and a scalar-valued
fluid pressure, i.e., $\{\hv,\hu_f,\hu_s,\hw, \hp_f\} 
\in \{ \hv^D + \hV^0\}
\times \{ \hu_f^D + \hV_{f,\hu}^0\}
\times \{ \hu_s^D + \hV_s^0 \}
\times \hV
\times\hL_f^0$, 
such that $\hv (0) = \hv^0$, $\hu_f (0) = \hu_f^0$, and 
$\hu_s(0) = \hu_s^0$ are satisfied, and for almost all
times $t\in\IT$ holds:
  \begin{eqnarray*}
    \begin{aligned}
&\text{Fluid/solid momentum}
\begin{cases}
     (\hJ\hrho_f  \partial_t \hv,\hpsi^v)_{\hOmega_f}  
      +(\hrho_f \hJ  (\hF^{-1}(\hv-\hw)\cdot\hnabla) \hv),
      \hpsi^v)_{\hOmega_f} 
      + (\hJ\hsigma_f\hF^{-T},\hnabla\hpsi^v)_{\hOmega_f}\\      
      + \langle \hrho_f \nu_f \hJ(\hF^{-T}\hnabla\hv^T\hn_f)\hF^{-T}, 
         \hpsi^v \rangle_{\hGamma_{\text{out}}}\\
      + (\hrho_s \partial_t \hv,\hpsi^v)_{\hOmega_s}  
      + (\hF\hSigma ,\hnabla\hpsi^v)_{\hOmega_s}
       =0\quad\forall\,\hpsi^v\in \hV^0, \\ 
\end{cases}\\
&\text{Fluid mesh motion (biharmonic/split)}
\begin{cases}   
      (\widehat\alpha \hnabla \hw|_{\hOmega_f},\hnabla\hpsi^u)_{\hOmega_f}
&=0\quad\forall\,\hpsi_f^u\in\hV_{f,\hu,\hGamma_i}^0 ,\\  
(\widehat\alpha\hw, \hpsi^w)_{\hOmega} 
- (\widehat\alpha\hnabla\hu_{f,s}, \hnabla\hpsi^w)_{\hOmega} 
& = 0\quad\forall\,\hpsi^w\in \hV
\end{cases}\\
&\text{Solid momentum, 2nd eq.}
\begin{cases}  
      \hrho_s (\partial_t\hu_s-\hv|_{\hOmega_s},\hpsi_s^u)_{\hOmega_s}
      & =0\quad\forall\,\hpsi_s^u\in \hL_s , \\ 
\end{cases}\\
&\text{Fluid mass conservation}
\begin{cases}  
      (\hdiv\,(\hJ\hF^{-1}\hv),\hpsi_f^p)_{\hOmega_f}
       & =0\quad\forall\,\hpsi_f^p\in \hL_f^0.
\end{cases}\\
    \end{aligned}
  \end{eqnarray*}  

The Neumann coupling 
conditions on $\hGamma_i$ are fulfilled in a variational way 
and cancel in monolithic modeling due to the global 
test space $\hV^0$ in which the test functions from 
both the fluid and the solid subdomains coincide on the interface. 
Thus, the condition
\begin{equation}
\langle \hJ\hsigma_f \hF^{-T}\hn_{f}, \hpsi^v \rangle_{\hGamma_i} 
+ \langle \hF\hSigma  \hn_{s},\hpsi^v 
\rangle_{\hGamma_i} = 0 
\quad\forall\,\hpsi^v \in \hV^0
\end{equation}
is implicitly contained in the above system.  
\end{form}

\section{Discretization}
\label{sec_dis}
In this section, we discuss temporal and spatial discretization 
of the forward problem. Our derivation contains many details 
on all terms of the FSI forward problem. 
The overall problem can be posed, however, in 
an abstract fashion, which facilitates the derivation 
of the backward-in-time adjoint problem in Section~\ref{sec_gradient}.

\subsection{Temporal discretization}
Our goal is to apply A-stable finite differences in time. 
Specifically, 
time discretization is based on a One-step-$\theta$ scheme 
as presented for the pure FSI problem, Formulation~\ref{eq:fsi:ale:harmonic},
in~\cite{Wi11}. 

In more detail, semi-discretization in time yields a sequence of generalized
steady-state problems that are completed by appropriate 
boundary values at every time step. Let 
\begin{equation*}
\IT = \{0\} \cup I_1 \cup \ldots \cup I_N
\end{equation*}
be a partition of the time interval $\IT$ into half open
subintervals $I_n :=(t_{n-1},t_n ]$ of (time step) size 
$k:= k_n := t_n - t_{n-1}$ with 
\begin{equation*}
0 = t_0 < \cdots < t_N = T.
\end{equation*}
Time derivatives are discretized with a backward 
difference quotient such that
\[
\partial_t
\hu \;\approx\;
\frac{\hu - \hu^{n-1}}{k}, \quad
\partial_t
\hv \;\approx\;
\frac{\hv - \hv^{n-1}}{k},
\]
where $\hu:=\hu^n:=\hu(t_n),\hv:=\hv^n:=\hv(t_n),
\hu^{n-1}:=\hu(t_{n-1}),\hv^{n-1}:=\hv(t_{n-1})$.
Furthermore, the mesh velocity $\partial_t \hcalA  = \hw$ is
numerically realized as $\hw = k^{-1}(\hu_f - \hu_f^{n-1})$.

\begin{form}[The time-discretized abstract problem]
\label{form_dis_abstract_problem_time}
We aim to
find $\hU^n=\{\hv^n, \hu_{f}^n, \hu_{s}^n,\hw^n,\hp_{f}^n\} \in \hX^0_{D}$, 
where $\hX_{D}^0 :=\{ \hv^D + \hV^0\} 
\times \{ \hu_{f}^D + \hV_{f,\hu}^0\} \times \hV_{s}^0
\times \hV
\times \hL_{f}^0$
and $\hX = \hV^0  \times \hV_{f,\hu,\hGamma_i}^0 \times 
\hV_{s}^0 \times \hV \times \hL_{f}^0$, 
for all $n=1,2,\ldots , N$ such that 
\begin{equation}
\label{eq_time_discretized_semilinear_form}
\hA(\hU^n)(\hPsi ) = 
0 
\quad\forall\, \hPsi \in \hX ,
\end{equation}
where the semi-linear form $\hA(\cdot)(\cdot)$ is split into
\begin{equation*}
\hA(\hU^n)(\hPsi ) := \hA_T(\hU^{n})(\hPsi ) + \hA_I(\hU^n)(\hPsi ) + 
\hA_E(\hU^n)(\hPsi ) + \hA_P(\hU^n)(\hPsi ).
\end{equation*} 
Details of this decomposition are provided in Definition~\ref{def_arrangement_semilinear_form}.
\end{form}

\begin{Definition}[Arranging the semi-linear form $\hA(\hU^n)(\hPsi )$ 
into groups]
\label{def_arrangement_semilinear_form}
We formally split the semilinear form into 
four categories: time equation terms (including the time
derivatives);
implicit terms (such as the fluid incompressibility and also 
the biharmonic mesh motion); pressure terms;
and finally all `standard' terms (e.g., stress terms, fluid convection).
We then obtain the decomposition:
\begin{equation}\label{Operator_Time_Equations_terms}
\begin{aligned}
  \hA_T (\hU)(\hPsi) 
  &= (\hJ \hrho_f \partial_t \hv ,\hpsi^v)_{\hOmega_f} 
  -(\hrho_f \hJ  (\hF^{-1}\hw\cdot\hnabla) \hv), \hpsi^v)_{\hOmega_f} 
   + (\hrho_s \partial_t \hv,\hpsi^v)_{\hOmega_s} 
  + (\hrho_s \partial_t\hu_s,\hpsi_s^u)_{\hOmega_s}, \\
  \hA_I (\hU)(\hPsi) 
  &=       (\widehat\alpha \hnabla \hw|_{\hOmega_f},\hnabla\hpsi^u)_{\hOmega_f}
+ (\widehat\alpha\hw, \hpsi^w)_{\hOmega} 
- (\widehat\alpha\hnabla\hu_{f,s}, \hnabla\hpsi^w)_{\hOmega} 
   + (\hdiv \,(\hJ\hF^{-1}\hv),\hpsi_f^p)_{\hOmega_f},  \\
 \hA_{P} (\hU)(\hPsi) &= (\hJ\hsigma_{f,p}
 \hF^{-T},\hnabla\hpsi^v)_{\hOmega_f}, \\
  \hA_E (\hU)(\hPsi) &= (\hrho_f \hJ  (\hF^{-1}\hv \cdot\hnabla) \hv),
  \hpsi^v)_{\hOmega_f} 
  + (\hJ\hsigma_{f,vu}\hF^{-T},\hnabla\hpsi^v)_{\hOmega_f}  \\
&\quad+ \langle \rho_f \nu \hJ (\hF^{-T}\hnabla v_f^T)\hF^{-T} \hn, \hpsi^v\rangle_{\hGamma_{\text{out}}} +(\hF\hSigma,\hnabla\hpsi^v)_{\hOmega_s}      
  - (\hrho_s \hv,\hpsi_s^u)_{\hOmega_s},
\end{aligned}
\end{equation}
where the fluid stress tensor $\hsigma_f$ is further split into $\hsigma_{f,vu}$, $\hsigma_{f,p}$:
\begin{align*}
\hsigma_{f,p} &= -\hp_f \hI, \quad
\hsigma_{f,vu} = \rho_f \nu_f (\hnabla\hv \hF^{-1} +
\hF^{-T}\hnabla\hv^T).
\end{align*}    
The (nonlinear) time derivative in  $\hA_T (\hU)(\hPsi)$ is approximated
by a backward difference quotient. For the time step 
$t_n \in \IT$, for $n=1,2,\ldots , N$ $(\, N\in\mathbb{N})$,
we compute $\hv := \hv^n , \hu_i := \hu_i^n \, (i=f,s)$ via
\begin{equation}
  \label{dis_backward_difference_time_derivative}
\begin{aligned}
  \hA_T (\hU^{n})(\hPsi) &\approx \frac{1}{k}\bigl(\hrho_f \hJ^{n,\theta} 
  (\hv - \hv^{n-1}),\hpsi^v \bigr)_{\hOmega_f} 
  -\frac{1}{k}\bigl( \hrho_f  (\hJ  \hF^{-1}
  (\hu_f - \hu_f^{n-1})\cdot\hnabla) \hv,
  \hpsi^v \bigr)_{\hOmega_f} \nonumber \\ 
  &\quad+ \frac{1}{k}\bigl(\hrho_s (\hv 
  - \hv^{n-1}),\hpsi^v\bigr)_{\hOmega_s} 
  + \frac{1}{k}\bigl(\hrho_s (\hu_s - \hu_s^{n-1}),\hpsi^u\bigr)_{\hOmega_s}\\
  &=: \frac{1}{k} \hA_{T,k}(\hU^{n},\hU^{n-1,t},\hPsi),
\end{aligned}
\end{equation}
where we introduce the  parameter $\theta\in [0,1]$. Furthermore, we use
\begin{equation*}
 \hJ^{n,\theta} = \theta\hJ^n + (1-\theta)\hJ^{n-1}, 
\end{equation*}
and $\hu_i^n :=\hu_i(t_n )$, $\hv^n :=\hv(t_n )$, 
and $\hJ:= \hJ^n := \hJ(t_n )$.
In our computations in Section~\ref{sec_tests}, we always consider $\hJ^{n,0.5}$.
The former time step is given 
by $\hv^{n-1}$, etc. for $i=f,s$.
\end{Definition}

\begin{form}
\label{form_time_discretized}
Let the previous time step solution 
$\hU^{n-1}=\{\hv^{n-1}, \hu_f^{n-1}, 
\hu_s^{n-1},\hw^{n-1},\hp_f^{n-1}\}$ and the 
time step $k:=k_n = t_n - t_{n-1}$ be given. 
In order to solve~\eqref{eq_time_discretized_semilinear_form}, we seek
$\hU^{n}=\{\hv^{n}, \hu_f^{n}, 
\hu_s^{n},\hw^n,\hp_f^{n}\}\in \hX^0$ by employing one-step-$\theta$ splitting:
\begin{align}
\label{computable_one_step_theta_scheme}
 \hA_{T,k} (\hU^{n},\hU^{n-1})(\hPsi) + \theta k\hA_E (\hU^n)(\hPsi) 
+ k \hA_P (\hU^n)(\hPsi) 
+ k \hA_I (\hU^n)(\hPsi) \;=\;  &- (1-\theta) k \hA_E (\hU^{n-1})(\hPsi).
\end{align}
The concrete scheme depends on the choice of the 
parameter $\theta\in [0,1]$. For $\theta = 1$ we obtain the strongly A-stable
backward Euler scheme (BE). 
If $k < 0.5$, for $\theta=0.5+k$, we obtain the 2nd order (shown for linear parabolic problems
in~\cite{Ra86, HeRa90}), A-stable, globally stabilized, 
Crank-Nicolson scheme (CNs).
\end{form}

\begin{Remark}
Formulation~\ref{form_time_discretized} is still nonlinear and continuous on the spatial
level. 
\end{Remark}

\subsection{Spatial discretization}
The time discretized formulation 
is the starting point for the Galerkin 
discretization in space. To this end, we construct a finite dimensional
subspace $\hX_h^0 \subset \hX^0$ to find an approximate solution to 
the continuous problem.
As previously explained, in the context of our variational-monolithic formulation, the computations
are done on the reference configuration $\hOmega$.
We use two dimensional shape-regular
meshes. A mesh
consists of quadrilateral cells $\widehat K$. They perform
a non-overlapping cover of the computation domain $\widehat\Omega
\subset\mathbb{R}^d$, $d=2$. 
The corresponding mesh is given by ${\cal \widehat T}_h = \{ \widehat K \}$. 
The discretization parameter in the reference 
configuration is denoted by $\widehat h$ and is a cell-wise constant
that is given by the diameter $\widehat h_{\widehat K}$ of the cell $\widehat K$.

On ${\cal \widehat T}_h$, the conforming finite element space for 
$\{\hv_{h}, \hu_{f,h},\hu_{s,h}, \hp_{f,h}, \hw_{h}\}$ is denoted by 
the space $\hX_h\subset \hX$. 
For Navier-Stokes flow, i.e., $\{\hv_{h},\hp_{f,h}\}$, we prefer the biquadratic, discontinuous-linear
$Q_2^c / P_1^{dc}$ element. For the specific 
definitions of the single elements, we refer the 
reader to~\cite{Cia87}. 
The property of the $Q_2^c / P_1^{dc}$ element is
continuity of the velocity values across different mesh cells~\cite{GiRa1986}. 
However, the pressure is defined by discontinuous test functions. 
Therefore, this element preserves
local mass conservation, is of low order, gains the \textit{inf-sup stability},
and is an optimal choice for both fluid problems and 
fluid-structure interaction
problems. The two displacement variables, namely 
$\hu_h,\hw_h$ are discretized with $Q_2^c$ elements.

%

In total, the discretized forward problem consists of 
\begin{form}
Given
$\hU^0\in\hX$ finding $\hU = (\hU^n_h)_{n=1}^N \in \hX^N_h$ solving 
\begin{equation}
  \label{eq:fsi_fullydiscrete}
  \begin{aligned}
    \sum_{n=1}^N \Biggl( &\hA_{T,k} (\hU^{n}_h,\hU^{n-1}_h)(\hPsi^n_h) + \theta k\hA_E (\hU^n_h)(\hPsi^n_h) 
    \\&+ k \hA_P (\hU^n_h)(\hPsi^n_h) 
    + k \hA_I (\hU^n_h)(\hPsi^n_h) + (1-\theta) k \hA_E
    (\hU^{n-1}_h)(\hPsi^n_h)
    \Biggr) = 0 \qquad \forall\, (\hPsi^n_h)_{n=1}^N\in \hX^N_h.
  \end{aligned}
\end{equation}
This abstract formulation serves as basis to derive the adjoint 
state in Section~\ref{sec_gradient}.
\end{form}

\section{Gradient computation}
\label{sec_gradient}
We are interested in identifying material parameters, e.g., $\mu$
in~\eqref{eq:Sigma}.
To this end, we denote by $q \in \R^p$, with $p\geq 1$, 
the collection of these
parameters, and will define suitable cost functionals $\mathcal
J(q,\hU)$ to be minimized. To highlight the dependence of the equation
on the parameters $q$, we add an additional $q$ argument to the form
$\hA_E$, e.g., we consider $\hA_E(q,\hU^n_h)(\hPsi^n_h)$ instead of
$\hA_E(\hU^n_h)(\hPsi^n_h)$ in~\eqref{eq:fsi_fullydiscrete}.

Assuming that~\eqref{eq:fsi_fullydiscrete} admits a unique solution
for any given $q$, we can obtain an optimality system 
by the standard Lagrange formalism, see,
e.g.,~\cite{Lions:1971,Troeltzsch:2010,HinzePinnauUlbrichUlbrich:2009}.
For a rigorous proof of the required differentiability properties,
some progress has been made for stationary FSI-problems
in~\cite{WiWo19}. A rigorous derivation of the corresponding adjoints
in the context of shape optimization can be found
in~\cite{HaubnerUlbrichUlbrich:2019}.

The formal Lagrange approach provides an adjoint equation
to~\eqref{eq:fsi_fullydiscrete} as
\begin{form}
Find $\hZ \in
\hX^N_h$ solving
\begin{equation}
  \label{eq:fsi_adjoint}
  \begin{aligned}
    \sum_{n=1}^N \Biggl( &\partial_{\hU^n}\hA_{T,k}
    (\hU^{n}_h,\hU^{n-1}_h)(\hPsi^n_h,\hZ^n_h)
    + \partial_{\hU^{n-1}}\hA_{T,k}
    (\hU^{n}_h,\hU^{n-1}_h)(\hPsi^{n-1}_h,\hZ^n_h)
    \\&+ \theta k\, \partial_{\hU^n}\hA_E (q,\hU^n_h)(\hPsi^n_h,\hZ^n_h) 
    + k\,\partial_{\hU^n} \hA_P (\hU^n_h)(\hPsi^n_h,\hZ^n_h) 
    + k\,\partial_{\hU^n} \hA_I (\hU^n_h)(\hPsi^n_h,\hZ^n_h) \\
    &+ (1-\theta) k \partial_{\hU^{n-1}}\hA_E
    (q,\hU^{n-1}_h)(\hPsi^{n-1}_h,\hZ^n_h)
    \Biggr) = \partial_{\widehat{U}}\mathcal J(q,\hU)(\hPsi) \qquad \forall\, (\hPsi^n_h)_{n=1}^N\in \hX^N_h.
  \end{aligned}
\end{equation}
Here $\partial_{\hU^n}\hA$ denotes the directional derivative of the
form $\hA$ with respect to its $\hU^n$ argument, and the first
argument of the second parentheses denotes the respective direction.
\end{form}

With this adjoint, we obtain the total derivative of the cost
functional $q \mapsto \mathcal J(q) :=\mathcal J(q,\hU)$ in a direction $\delta q$ as
\[
  \frac{\mathrm{d}}{\mathrm{d}q} J(q,\hU)\delta q = \partial_q
  \mathcal J(q,\hU)(\delta q) + \sum_{n=1}^N
  \Biggl(\theta k \,\partial_{q}\hA_E (q,\hU^n_h)(\delta q,\hZ^n_h)
  + (1-\theta) k\, \partial_{q}\hA_E(q,\hU^{n-1}_h)(\delta
  q,\hZ^n_h)\Biggr)
\]
allowing the calculation of the reduced gradient $\nabla J(q)
\in \R^p$ of the cost functional by
\begin{equation}\label{eq:gradient}
  (\nabla J(q),\delta q) = \frac{\mathrm{d}}{\mathrm{d}q}
  J(q,\hU)\delta q \qquad \forall\, \delta q \in \R^p
\end{equation}
cf., e.g.,~\cite{BeckerMeidnerVexler:2007}.

\section{Solution algorithms}
\label{sec_algo}
In order to minimize the cost functional $\mathcal J(q)$, 
we employ a standard globalized gradient method, i.e.,
\begin{algo}[Gradient method]\label{algo:gradient}
  Let $q^0\in \R^p$  be an initial guess, and pick parameters $\gamma \in
  (0,1/2)$ and $\beta \in (0,1)$. For $k=0,1,\dots$ until $\|\nabla_q
  \mathcal J(q^k)\|_Q < TOL$ iterate
  \begin{enumerate}
  \item Solve the 
  (nonlinear) primal problem~\eqref{eq:fsi_fullydiscrete} to
  obtain $\hU_h\in \hX^N_h$ using
  Algorithm~\ref{algo_residual_based_Newton}. 
  \item Solve the (linear) adjoint problem~\eqref{eq:fsi_adjoint} to obtain
    $\hZ_h\in \hX^N_h$.
  \item Compute the gradient  $\nabla \mathcal J(q^{k})$
    using~\eqref{eq:gradient}.
  \item Find the largest $l \in \{0,1,\ldots\}$ such that (Armijo-rule) 
    \[
      \mathcal J(q^{k} - \beta^l \nabla \mathcal J(q^{k})) \le \mathcal J(q^{k}) - \gamma \beta^l \|\nabla \mathcal J(q^{k})\|^2
    \]
    holds and set $\beta_k = \beta^l$.
  \item Update
    \[
    q^{k+1} = q^{k} - \beta_k \nabla \mathcal J(q^{k}).
    \]
  \end{enumerate}
\end{algo}


In the first step of the previous algorithm for solving the nonlinear 
primal problem, we employ 
the following Newton method. At each time point
the following problem is given:
\begin{equation*}
\widehat A (\widehat U_h^{n})(\widehat\Psi ) = 
0 \quad \forall\, \widehat\Psi\in \widehat X_h.
\end{equation*}

\begin{algo}[Residual-based Newton's method]
\label{algo_residual_based_Newton}
We omit $h$ and $n$ for the convenience of the reader.
Choose an initial Newton guess $\hU^0\in \hX$. For the iteration steps 
$k=0,1,2,3,\ldots$:
\begin{enumerate}
\item Find $\delta \hU^k\in \hX^0$ such that
\begin{align}
\hA'(\hU^k)(\delta \hU^k, \hPsi) &= 
- \hA(\hU^k)(\hPsi)\quad\forall\,\hPsi\in
\hX, \label{linear_system_in_Newton}  \\
\hU^{k+1} &= \hU^k + \lambda_k \delta \hU^k, \label{linear_system_in_Newton_sol}
\end{align}
for $\lambda_k = 1$. The arising linear equations are solved 
with a direct method; namely UMFPACK~\cite{DaDu97}. This choice is justified 
since the spatial numbers of degrees of freedom is moderate in our numerical 
examples. Moreover, in order to save computational cost, we adopt 
simplified Newton steps; i.e., the matrix $\hA'(\hU^k)(\delta \hU^k, \hPsi)$ is only 
rebuild when $\lambda_k^l < 1$ (defined below) or 
 $\|\hA(\hU^k)\| \in [0.001,1] \|\hA(\hU^{k-1})\|$.
\item The criterion for convergence is the contraction of the residuals:
\begin{equation}
\label{eq_A_k1_Ak}
\| \hA(\hU^{k+1})\| < \|\hA(\hU^k)\|.
\end{equation}
\item If~\eqref{eq_A_k1_Ak} is violated, re-compute in~\eqref{linear_system_in_Newton_sol}
$\hU^{k+1}$ by choosing $\lambda_k = 0.6^l$, and
compute for $l=1,...,l_M$ (e.g. $l_M = 5$) a new solution
\[
\hU^{k+1} = \hU^k + \lambda_k^l \delta \hU^k
\]
until~\eqref{eq_A_k1_Ak} is fulfilled for a $l^*<l_M$ or $l_M$ is reached.
In the latter case, no convergence is obtained and the program aborts.
\item In case of $l^*<l_M$ we check next the (relative) stopping criterion: 
\[
\|\hA(\hU^{k+1})\| \leq \|\hA(\hU^{0})\|\, TOL_N. 
\]
If this is criterion is fulfilled, set $\hU^n := \hU^{k+1}$.
Otherwise, we increment $k\to k+1$ and goto Step~1.
\end{enumerate}
\end{algo}

\section{Numerical tests}
\label{sec_tests}
We conduct several numerical tests in this section. These
are implemented in the open-source 
package DOpElib~\cite{dope,DOpElib} 
using the finite elements of deal.II~\cite{dealII90}.

An open-source implementation of~\eqref{eq:fsi_fullydiscrete}, 
used as basis for
our computations, can be found in DOpElib~\cite{dope,DOpElib} in
\texttt{Examples/PDE/InstatPDE/Example2}. 

\subsection{Material parameters}
The material parameters for the FSI-1 and FSI-3 tests are chosen as proposed
in~\cite{HrTu06b,BuSc10} 
and listed in Table~\ref{material_param}. 
The parameters for the flapping example are a mixture of 
\cite{GiCaBoHa10,Wi12_fsi_eale_heart,FaiWi18}.

\renewcommand{\arraystretch}{1.5}

\begin{table}[h!]
\centering
 \begin{tabular}{ l | r | r | r }
		& FSI-1 			& FSI-3 							& Flapping 		\\ \hline
$\nu_f$ 	& $10^{-3} \si{\frac{m^2}{s}}$& $10^{-3}\si{\frac{m^2}{s}}$	& $10^{-1} \si{\frac{cm^2}{s}}$\\	
$\mu_s$		& $0.5\cdot 10^6 \si{\frac{kg}{m\, s^2}}$	& $2.0\cdot10^6 \si{\frac{kg}{m\, s^2}}$& $1.0\cdot 10^9 \si{\frac{g}{cm\, s^2}}$\\
$\nu_s$         & $0.4$ & $0.4$ & $0.4$ \\
$\hrho_s$	& $10^3 \si{\frac{kg}{m^3}}$	& $10^3 \si{\frac{kg}{m^3}}$	& $10^2 \si{\frac{g}{cm^3}}$	\\	
$\hrho_f$ 	& $10^3 \si{\frac{kg}{m^3}}$	& $10^3 \si{\frac{kg}{m^3}}$   & $10^2 \si{\frac{g}{cm^3}}$\\ \hline
\end{tabular}
\caption{Material parameters for all test cases.
}\label{material_param}
\end{table}

\renewcommand{\arraystretch}{1.0}

\subsection{Example 1: Parameter estimation within the FSI~1 benchmark}
In this first numerical example, we consider a quasi-stationary 
setting based on the FSI 1 benchmark~\cite{HrTu06b,BuSc10}.
The optimization problem reads: identify 
the Lam\'e parameter $\mu_s$ from measurements 
of the beam-tip displacement at $A:=(0.6,0.2)$. The exact values $u_d$
are taken from a reference solution with $\mu = 0.5\cdot 10^6 \si{\frac{kg}{m\, s^2}}$.

The forward problem is solved with the backward Euler scheme, i.e.,
$\theta = 1$, since the configuration is stationary and we only use a
time-dependent method to find the stationary limit.

\subsubsection{Cost functional}
The cost functional reads:
\[
J(q,\hU) = \frac{1}{2}(\hu_1(A,T)-u_{d})^2 +
  \frac{\alpha}{2}|q-q_d|^2,
\]
with $q_d = 2.27007 \cdot 10^{-5}$.

\subsubsection{Configuration}
\label{sec_config_fsi_1}
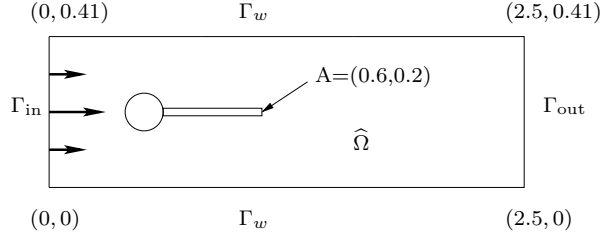
\begin{figure}[h!]
\centering
\begin{picture}(0,0)%
\includegraphics{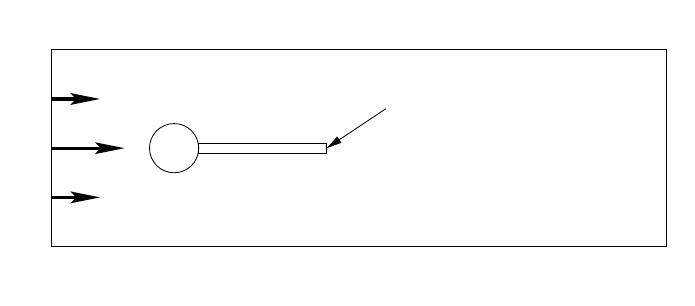}%
\end{picture}%
\setlength{\unitlength}{2072sp}%
\begingroup\makeatletter\ifx\SetFigFont\undefined%
\gdef\SetFigFont#1#2{%
  \fontsize{#1}{#2pt}%
  \selectfont}%
\fi\endgroup%
\begin{picture}(6330,2785)(3136,-4652)
\put(9001,-4561){\makebox(0,0)[lb]{\smash{{\SetFigFont{8}{9.6}{\color[rgb]{0,0,0}$(2.5,0)$}%
}}}}
\put(9001,-2086){\makebox(0,0)[lb]{\smash{{\SetFigFont{8}{9.6}{\color[rgb]{0,0,0}$(2.5,0.41)$}%
}}}}
\put(3376,-2086){\makebox(0,0)[lb]{\smash{{\SetFigFont{8}{9.6}{\color[rgb]{0,0,0}$(0,0.41)$}%
}}}}
\put(3376,-4561){\makebox(0,0)[lb]{\smash{{\SetFigFont{8}{9.6}{\color[rgb]{0,0,0}$(0,0)$}%
}}}}
\put(6751,-2851){\makebox(0,0)[lb]{\smash{{\SetFigFont{8}{9.6}{\color[rgb]{0,0,0}A=(0.6,0.2)}%
}}}}
\put(7201,-3661){\makebox(0,0)[lb]{\smash{{\SetFigFont{8}{9.6}{\color[rgb]{0,0,0}$\widehat{\Omega}$}%
}}}}
\put(5851,-2086){\makebox(0,0)[lb]{\smash{{\SetFigFont{8}{9.6}{\color[rgb]{0,0,0}$\Gamma_w$}%
}}}}
\put(5851,-4561){\makebox(0,0)[lb]{\smash{{\SetFigFont{8}{9.6}{\color[rgb]{0,0,0}$\Gamma_w$}%
}}}}
\put(3151,-3211){\makebox(0,0)[lb]{\smash{{\SetFigFont{8}{9.6}{\color[rgb]{0,0,0}$\Gamma_{\text{in}}$}%
}}}}
\put(9451,-3211){\makebox(0,0)[lb]{\smash{{\SetFigFont{8}{9.6}{\color[rgb]{0,0,0}$\Gamma_{\text{out}}$}%
}}}}
\end{picture}%

\caption{FSI-1 and FSI-3 benchmarks (Examples 1 and 2): flow around cylinder with elastic beam with circle-center $C=(0.2,0.2)$ and radius $r=0.05$.}
\label{configuration_csm_and_fsi_2D}
\end{figure}

The geometry of the FSI-1 and FSI-3 settings
are displayed in Figure~\ref{configuration_csm_and_fsi_2D}.
An elastic beam is attached to a cylinder and
is surrounded by an incompressible fluid. The initial geometry is once 
uniformly refined in space.

On the
cylinder and outer boundary $\Gamma_w$ we enforce zero Dirichlet boundary
conditions for $\hv$ and $\hu$. On the outflow boundary $\Gamma_{\text{out}} $ we prescribe the do-nothing
outflow condition~\cite{HeRaTu96}. The inflow profile on $\Gamma_{\text{in}}$ is given by:
\[
\hv(0,y):= 1.5\, y\, (0.41-y)\, \frac{4}{0.41^2}\, v_{\text{mean}}(t).
\]
The mean inflow $v_{\text{mean}}(t)$ is $0.2\, \si{m/s}$ for Example 1
(FSI 1) and $2.0\, \si{m/s}$ in Example 2 (FSI 3).
In the FSI 1 test case, we compute $n=25$ time steps using $k=1\, \si{s}$ and 
in the FSI 3 example, we work with $k=0.001\, \si{s}$ with $T=0.6\, \si{s}$ corresponding 
to $n=6000$ time steps.

\subsubsection{Discussion of the FSI 1 findings}
Our results for three different configurations are displayed 
in the Tables~\ref{Tab:FSI1-1},~\ref{Tab:FSI1-2}, and~\ref{Tab:FSI1-3}.
In the first run with $\alpha = 0.001$, the algorithm converges slowly in order 
to estimate $q^k$ and to reduce the cost functional 
$\mathcal J(q^k)$. The main reason is due to the low regularization, which 
is confirmed by two further runs with $\alpha = 0.1$ and $1$.

\begin{table}[h!]
   \caption{Optimization results for the FSI 1 example with $\alpha = 0.001$ and $q_d = \cdot 10^6$. The initial Residual in $q_0=5\,000$ is $|\nabla
\mathcal J(q^0)| = 987.5$}
\label{Tab:FSI1-1}
   \centering   
   \begin{tabular}{@{}llll@{}}
\hline
Iter     & $\mathcal J(q^k)$         & $ q^k          $ &   $\frac{|\nabla
\mathcal J(q^k)|}{|\nabla
\mathcal J(q^0)|}$ \\
0   & $4.913 \cdot 10^{8} $ & $5000   $ & $ 1.0000\cdot 10^{-0} $\\
1   & $4.9033\cdot 10^{8} $ & $5987.54$ & $ 9.9901\cdot 10^{-1} $\\
2   & $4.8936\cdot 10^{8} $ & $6974.11$ & $ 9.9802\cdot 10^{-1} $\\
3   & $4.8838\cdot 10^{8} $ & $7959.69$ & $ 9.9703\cdot 10^{-1} $\\
$\vdots$	 & $\vdots$    &  $\vdots$           & $\vdots$ \\
101 & $4.0201\cdot 10^{8} $ & $99950.5$ & $ 9.0457\cdot 10^{-1} $\\
102 & $4.0121\cdot 10^{8} $ & $100844 $ & $ 9.0367\cdot 10^{-1} $\\
103 & $4.0042\cdot 10^{8} $ & $101736 $ & $ 9.0278\cdot 10^{-1} $\\
$\vdots$	 & $\vdots$    &  $\vdots$           & $\vdots$ \\
198 & $3.3157\cdot 10^{8} $ & $182600 $ & $ 8.2151\cdot 10^{-1} $\\
199 & $3.3091\cdot 10^{8} $ & $183411 $ & $ 8.2069\cdot 10^{-1} $\\
200 & $3.3025\cdot 10^{8} $ & $184222 $ & $ 8.1988\cdot 10^{-1} $\\
$\vdots$	 & $\vdots$    &  $\vdots$           & $\vdots$ \\
\hline
\end{tabular}
\end{table}

In Table~\ref{Tab:FSI1-2}, the value of $\alpha$ is enlarged to $0.1$. 
Here, in $155$ gradient iterations, the cost functional is 
reduced by a order to $10^{14}$ from an initial control $q^0 = 5000$ 
to $q^{155} = 10^6$. 

\begin{table}[h!]
   \caption{Optimization results for the FSI 1 example with $\alpha = 0.1$ and $q_d = 10^6$. The initial Residual in $q_0=5\,000$ is $|\nabla
\mathcal J(q^0)| = 9.875 \cdot 10^4$}
\label{Tab:FSI1-2}
   \centering   
   \begin{tabular}{@{}llll@{}}
\hline
Iter     & $\mathcal J(q^k)$         & $ q^k          $ &   $\frac{|\nabla
\mathcal J(q^k)|}{|\nabla
\mathcal J(q^0)|}$ \\
0   & $ 4.913 \cdot 10^{10} $ & $5000   $ & $ 1.0000\cdot 10^{-0} $\\
1   & $ 3.9862\cdot 10^{10} $ & $103754 $ & $ 9.0075\cdot 10^{-1} $\\
2   & $ 3.2342\cdot 10^{10} $ & $192707 $ & $ 8.1135\cdot 10^{-1} $\\
3   & $ 2.6241\cdot 10^{10} $ & $272832 $ & $ 7.3082\cdot 10^{-1} $\\
$\vdots$	 & $\vdots$    &  $\vdots$           & $\vdots$ \\
101 & $ 3.3216\cdot 10^{ 1} $ & $999974 $ & $ 2.6001\cdot 10^{-5} $\\
102 & $ 2.6950\cdot 10^{ 1} $ & $999977 $ & $ 2.3421\cdot 10^{-5} $\\
103 & $ 2.1865\cdot 10^{ 1} $ & $999979 $ & $ 2.1096\cdot 10^{-5} $\\
$\vdots$	 & $\vdots$    &  $\vdots$           & $\vdots$ \\
154 & $ 5.1211\cdot 10^{-4} $ & $10^{6} $ & $ 1.0210\cdot 10^{-7} $\\
155 & $ 4.1550\cdot 10^{-4} $ & $10^{6} $ & $ 9.1962\cdot 10^{-8} $\\
\hline
\end{tabular}
\end{table}

Increasing further $\alpha$ to $1$ (Table~\ref{Tab:FSI1-2}) yields 
a reduction in $\mathcal J(q^k)$ from about $10^{11}$ to $10^{-6}$.
The gradient algorithm converges in $5$ iterations.

\begin{table}[h!]
   \caption{Optimization results for the FSI 1 example with $\alpha = 1$ and $q_d = 500\,000$. The initial Residual in $q_0=5\,000$ is $|\nabla
\mathcal J(q^0)| = 4.913\cdot 10^5$}
\label{Tab:FSI1-3}
   \centering   
   \begin{tabular}{@{}llll@{}}
\hline
Iter     & $\mathcal J(q^k)$         & $ q^k          $ &   $\frac{|\nabla
\mathcal J(q^k)|}{|\nabla
\mathcal J(q^0)|}$ \\
0 & $ 1.216 \cdot 10^{11} $ & $ 5000   $ & $ 1.0000\cdot 10^{-0} $ \\
1 & $ 6.8268\cdot 10^{ 6} $ & $ 496291 $ & $ 7.4929\cdot 10^{-3} $ \\
2 & $ 3.8328\cdot 10^{ 2} $ & $ 499972 $ & $ 5.6144\cdot 10^{-5} $ \\
3 & $ 2.1519\cdot 10^{-2} $ & $ 500000 $ & $ 4.2068\cdot 10^{-7} $ \\
4 & $ 1.2082\cdot 10^{-6} $ & $ 500000 $ & $ 3.1522\cdot 10^{-9} $ \\
\hline
\end{tabular}
\end{table}

\subsection{Example 2: Optimal control within the FSI~3 benchmark}
In this second numerical test, we employ 
the same geometry as in Example 1. The material 
parameters and boundary data can be found in 
Table~\ref{material_param} and Section~\ref{sec_config_fsi_1}.
We now consider an optimal control problem in which 
$\mu_s$ is detected in such a way to match the displacement value
at the beam tip at $(0.6,0.4)$ obtained by the FSI 1 simulation
in Example 1. Since this numerical test is nonstationary 
with periodic solutions in the original forward run, we use
the shifted Crank-Nicolson time-stepping scheme with minimal 
numerical dissipation.

\subsubsection{Cost functional}
The cost functional is given by:
\[
J(q,\hU) = \frac{1}{2}(\hu_1(A,T)-u_{d})^2 +
  \frac{\alpha}{2}|q-q_d|^2
\]
with $q_d = 2.27007 \cdot 10^{-5}$, i.e., the functional
is similar to Example 1 and the reference value 
comes from FSI 1, but not FSI 3.

\subsubsection{Discussion of the FSI 3 findings}
\begin{figure}[h!]
  \begin{center}
  \includegraphics[width = 8cm]{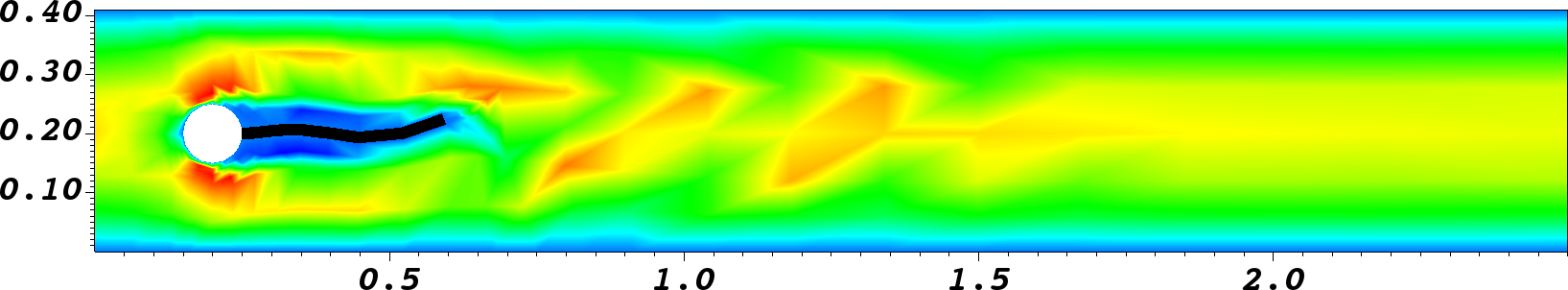}
  \includegraphics[width = 8cm]{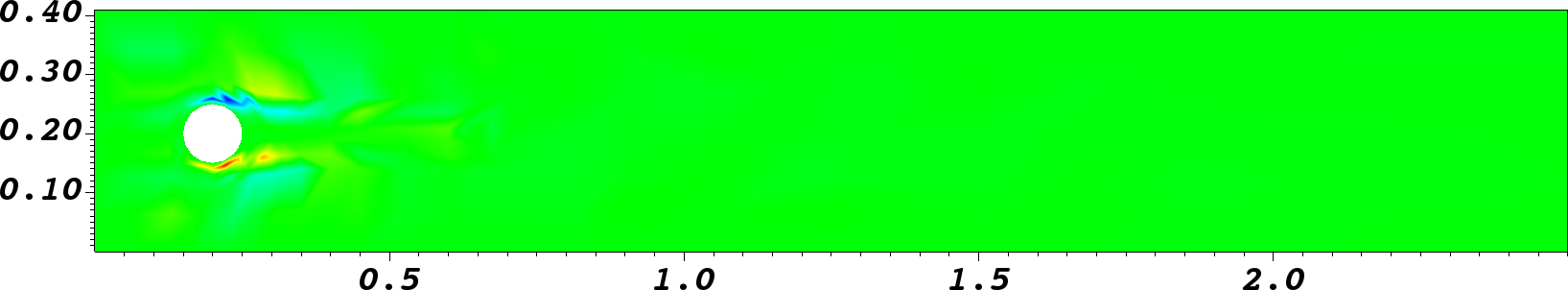}\\
  \includegraphics[width = 8cm]{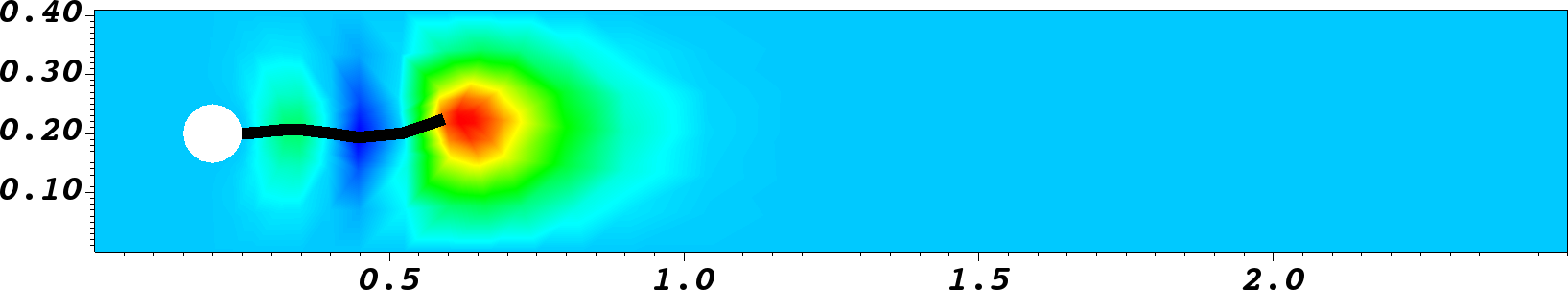}
  \includegraphics[width = 8cm]{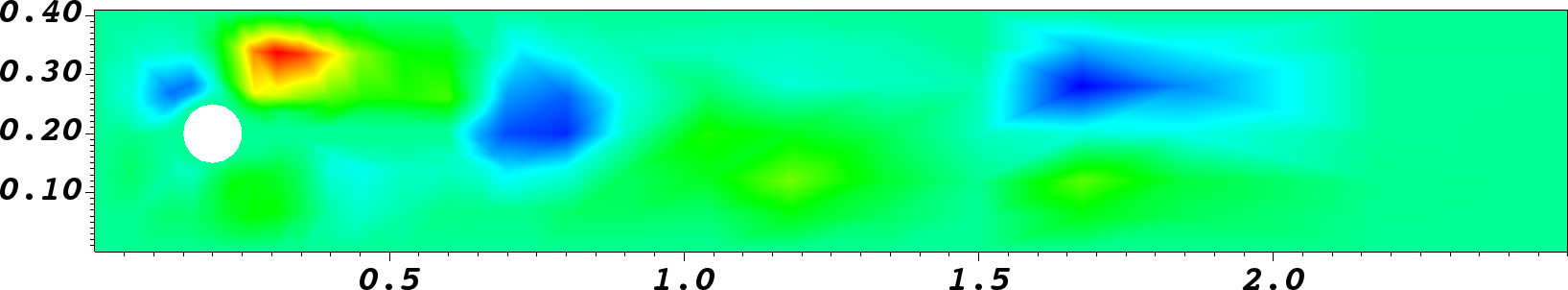}\\
  \includegraphics[width = 8cm]{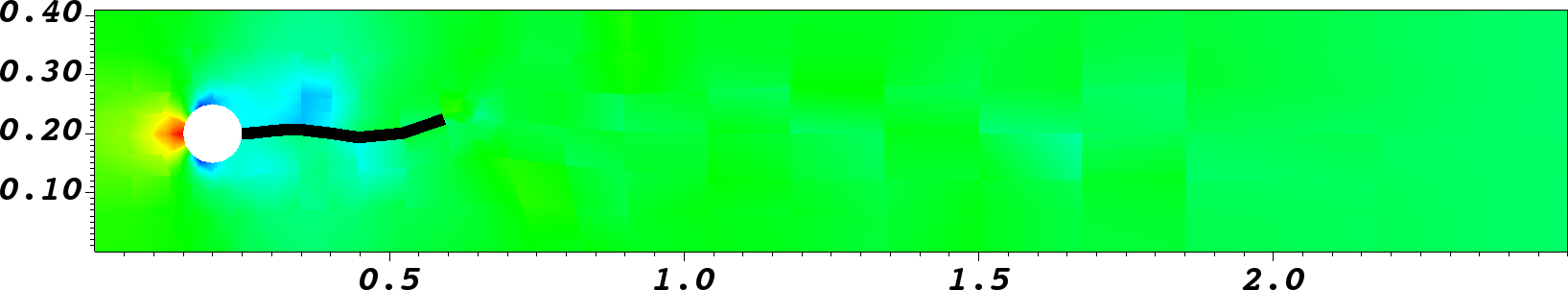}
  \includegraphics[width = 8cm]{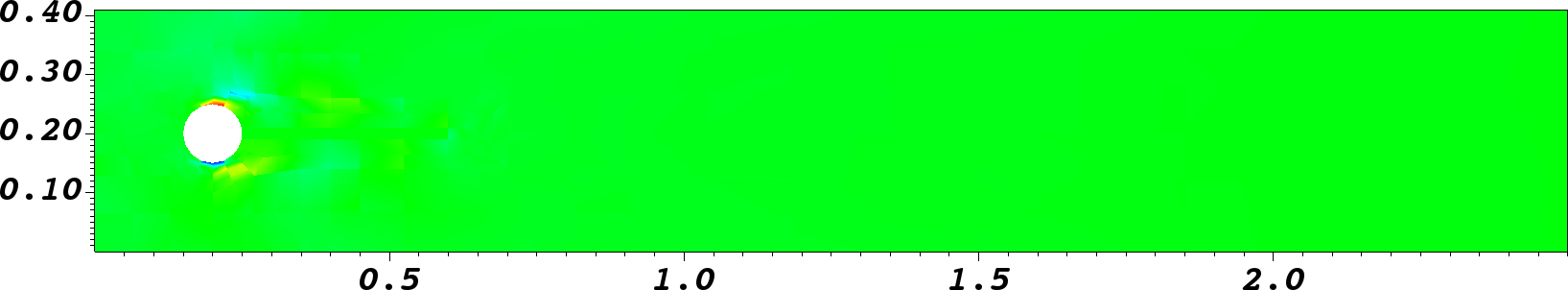}
    \caption{
      Example 2: At $T=5s$ (time step No. $5000$): $v_x(t), u_y(t)$ and $p(t)$ in the 
deformed configuration $\Omega(t)$. Left column: the primal states are shown. Right column: the corresponding 
adjoint states are shown.} 
    \label{example_2_plots}
  \end{center}
 \end{figure}

Graphical plots of the solution are provided in Figure 
\ref{example_2_plots}. Our quantitative results are shown
in Table~\ref{Tab:FSI3}. The gradient algorithm converges 
in $29$ iterations in which the cost 
functional is reduced by $10^3$ and the 
control is approximated by $q^{29} = 572\, 378$.

\begin{table}[h!]
   \caption{Optimization results for the FSI-3 example with $\alpha = 0.1$ and $q_d = 500\,000$. The initial Residual in $q_0=2\cdot 10^{6}$ is $1.489\cdot 10^5$}
\label{Tab:FSI3}
   \centering   
   \begin{tabular}{@{}llll@{}}
\hline
Iter     & $\mathcal J(q^k)$            & $ q^k          $ &    $\frac{|\nabla
\mathcal J(q^k)|}{|\nabla
\mathcal J(q^0)|}$\\
0	 & $1.117 \cdot 10^{11} $	& $ 2      \cdot 10^{6} $ & $ 1.0000\cdot 10^{-0} $ \\      	
1	 & $9.0593\cdot 10^{10} $	& $ 1.85112\cdot 10^{6} $ & $ 9.0075\cdot 10^{-1} $ \\    
2	 & $7.3502\cdot 10^{10} $	& $ 1.71702\cdot 10^{6} $ & $ 8.1135\cdot 10^{-1} $ \\    
3	 & $5.9636\cdot 10^{10} $	& $ 1.59623\cdot 10^{6} $ & $ 7.3082\cdot 10^{-1} $ \\    
4	 & $4.8386\cdot 10^{10} $	& $ 1.48743\cdot 10^{6} $ & $ 6.5829\cdot 10^{-1} $ \\    
$\vdots$	 & $\vdots$    &  $\vdots$           & $\vdots$ \\
28	 & $3.2041\cdot 10^{ 8} $	& $ 580353 $              & $ 5.3569\cdot 10^{-2} $ \\
29	 & $2.5996\cdot 10^{ 8} $	& $ 572378 $              & $ 4.8252\cdot 10^{-2} $ \\
\hline
\end{tabular}
\end{table}

\subsection{Example 3: Two-dimensional flapping membranes}
In this third example, we consider two-dimensional 
flap dynamics. This test is a challenge because of the thin flaps 
and the mesh regularity. The original setups for forward 
simulations were inspired by~\cite{GiCaBoHa10}. Our configuration 
here is a further extension, towards FSI-optimization, of~\cite{Wi12_fsi_eale_heart} 
and \cite{FaiWi18}.

\subsubsection{Cost functional}
The cost functional is given by:
\[
J(q,\hU) = F(\hGamma_{\text{opt}}, T) + \frac{\alpha}{2}|q-q_d|^2
\]
where $T$ is the end time value as in the other examples 
and $F(\cdot)$ is the drag functional defined as
\[
F(\hGamma_{\text{opt}}, T) := \int_{\hGamma_{\text{opt}}} (\hsigma_f \cdot \hn)\cdot e_1 
\, \mathrm{d}s
\]
where $\hn$ is the unit normal vector pointing outward of the domain $\hOmega_s$ 
and $e_1$ the first unit vector in $\mathbb{R}^2$.
The boundary part, where the drag is evaluated is 
\[
\hGamma_{\text{opt}} := \{ 2\leq x \leq 8; \; y=0 \}.
\]
Moreover, we notice that we only control $\mu$ in the
  valves, while in the rest of the solid, the value is as 
in Table \ref{material_param}.

\subsubsection{Configuration}
\begin{figure}[h!]
  \begin{center}
    \includegraphics[width = 10cm]{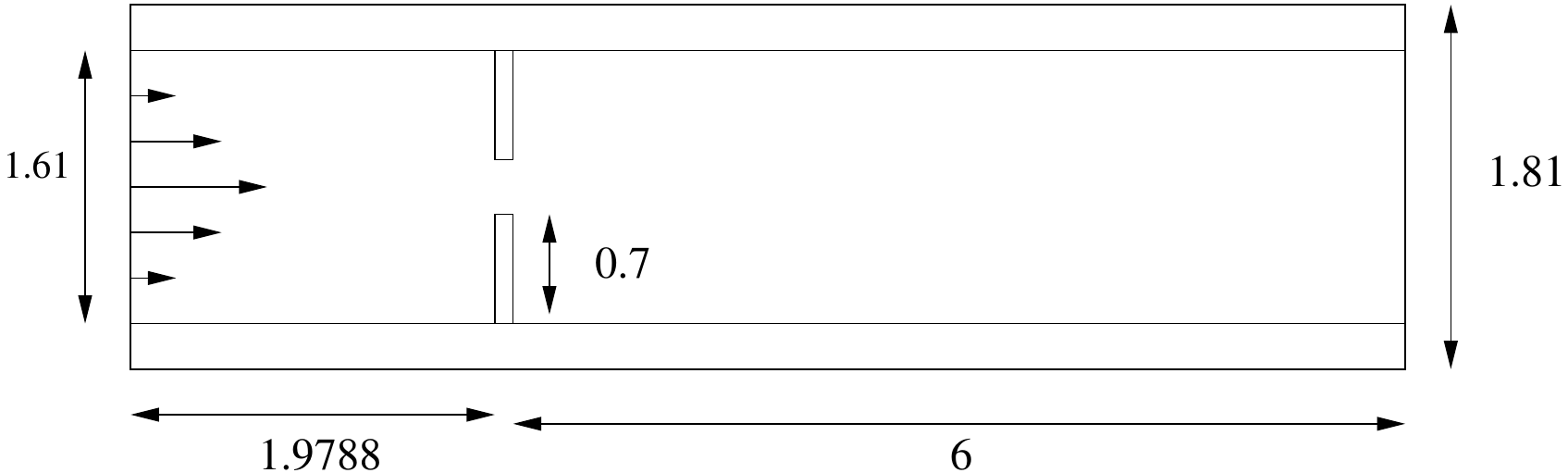}
    \caption{
      Example 3: Configuration. All data given in $\si{cm}$.} 
    \label{configuration_gil}
  \end{center}
 \end{figure}

\begin{figure}[h!]
  \begin{center}
    \includegraphics[width = 14cm]{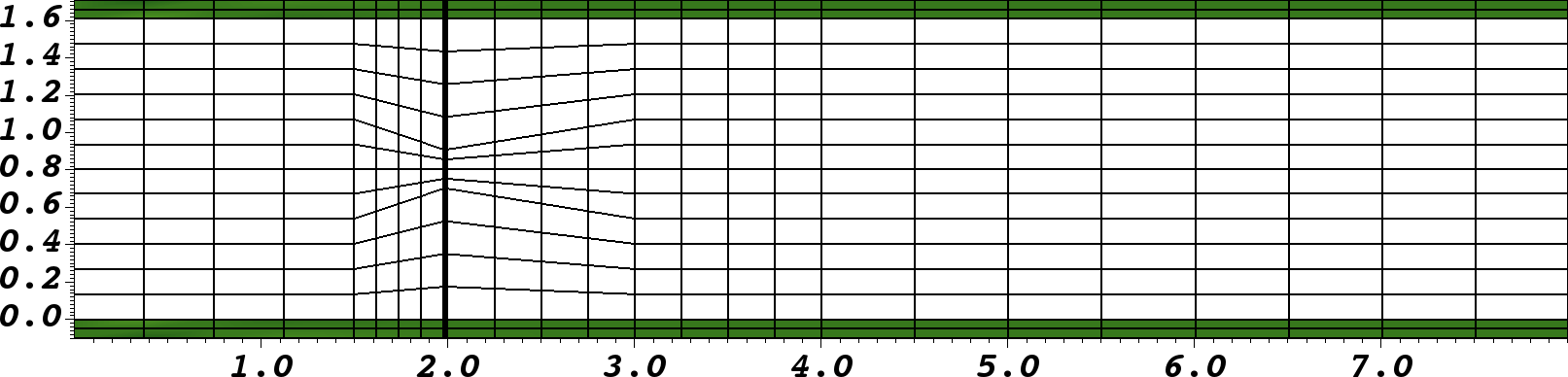}
    \caption{
      The mesh for the flapping membranes example at the initial time step.
      All geometric values are given in $\si{cm}$.
      The solid boundaries are colored in dark green. The flaps are located 
      at $1.9788\, \si{cm} \leq x \leq 2.0\, \si{cm}$.} 
    \label{initial mesh valves}
  \end{center}
 \end{figure}

The geometry is shown in Figure~\ref{configuration_gil}.
The initial mesh is once uniformly refined yielding 
the mesh shown in Figure~\ref{initial mesh valves}.

%

\begin{figure}[h!]
   \begin{center}
	\includegraphics{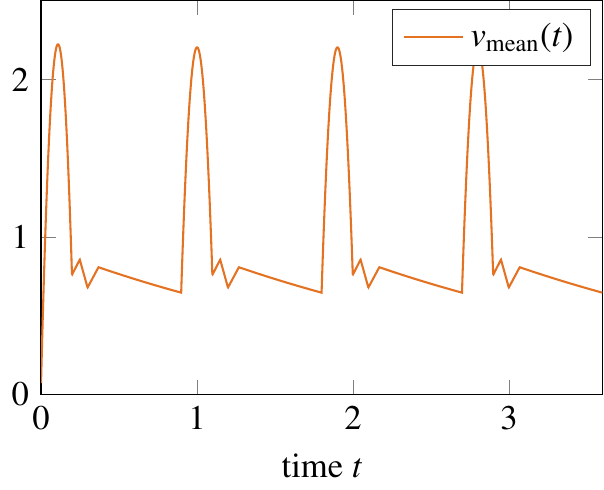}
 \caption{
Interpolated flow rate profile $\bar{v}(t)$ that 
is used to scale the inflow profile of the flapping membrane example.} 
 \label{velocity_profile_heart}
  \end{center}
 \end{figure}

On the inflow boundary, $\hGamma_{\text{in}} := \{ x=0; -0.1 \leq y \leq 1.61 \}$, we prescribe a parabolic inflow profile
\[
v(0,y):= 0.15 y (1.61-y) \frac{4}{1.61^2} v_{\text{mean}}(t) \quad\text{for } t\in I:=[0,0.9],
\]
where $v_{\text{mean}}(t)$ taken 
from Figure~\ref{velocity_profile_heart}.

At the outflow boundary the do-nothing outflow 
condition  $\hGamma_{\text{out}} $ is prescribed for $\hv$ and $\hp$, while 
the displacements are fixed there.
On the outer wall boundaries 
\[
  \hGamma_{\text{wall}} := \{ 0\leq x \leq 8; y=-0.1 \}\cup \{ 0\leq x
  \leq 8; y=-1.61 \}
\]
we use homogeneous Neumann conditions for the displacements and the velocity
in order to allow the solid to move freely.

The computations are performed on the time interval
$I= (0, 0.579375\,\si{s})$. 
The end time value $T=0.579375\,\si{s}$ is chosen such that the first 
maximal stress appears for the initial control $q^0$. For the
computations, the time
interval is split into $618$ time steps.

\subsubsection{Discussion of the flapping membrane findings}
\begin{figure}[h!]
  \begin{center}
  \includegraphics[width = 15cm]{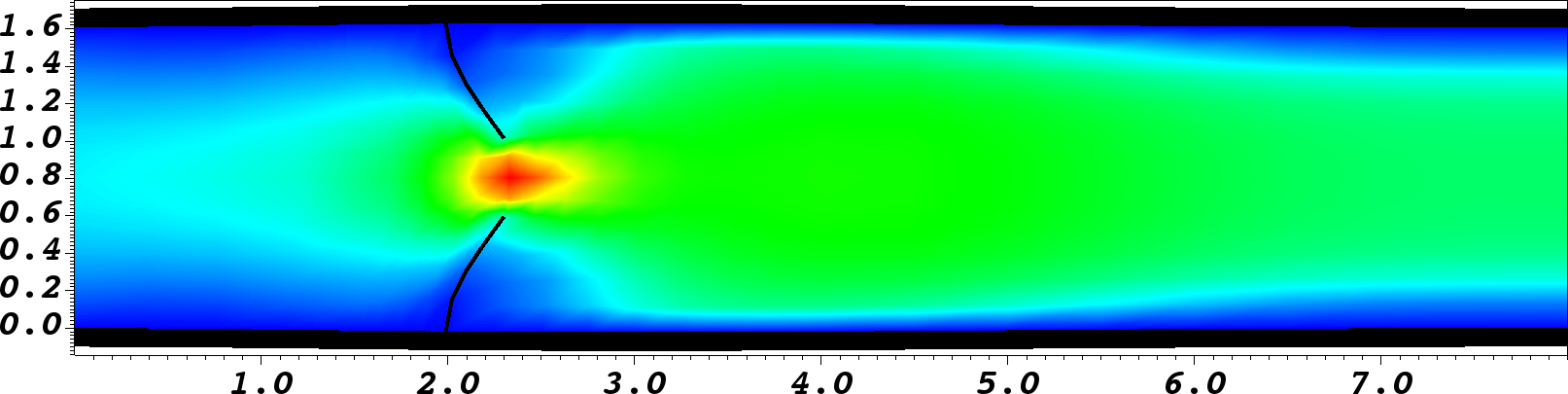}\\[1em]
  \includegraphics[width = 15cm]{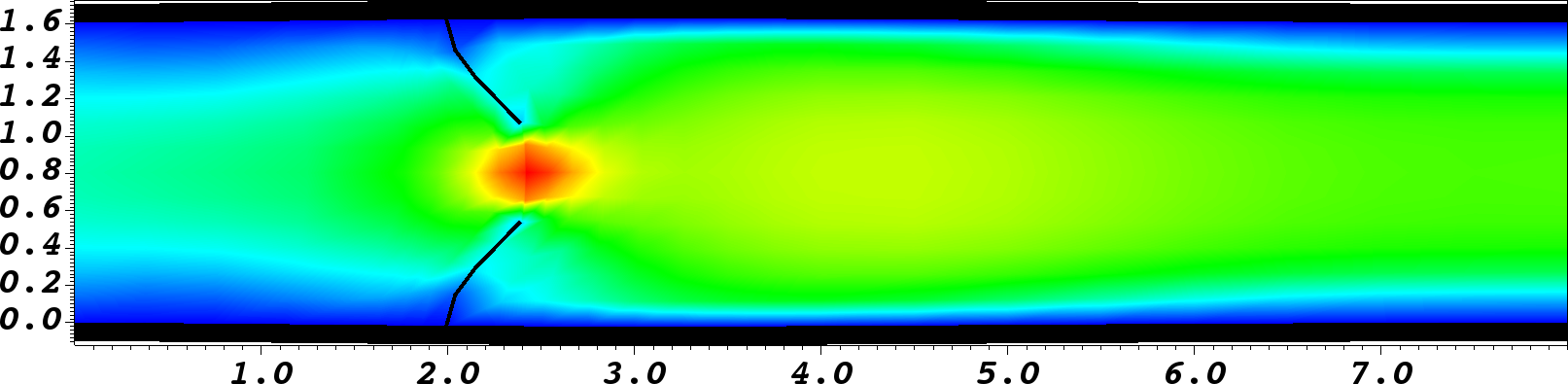}\\[1em]
  \includegraphics[width = 15cm]{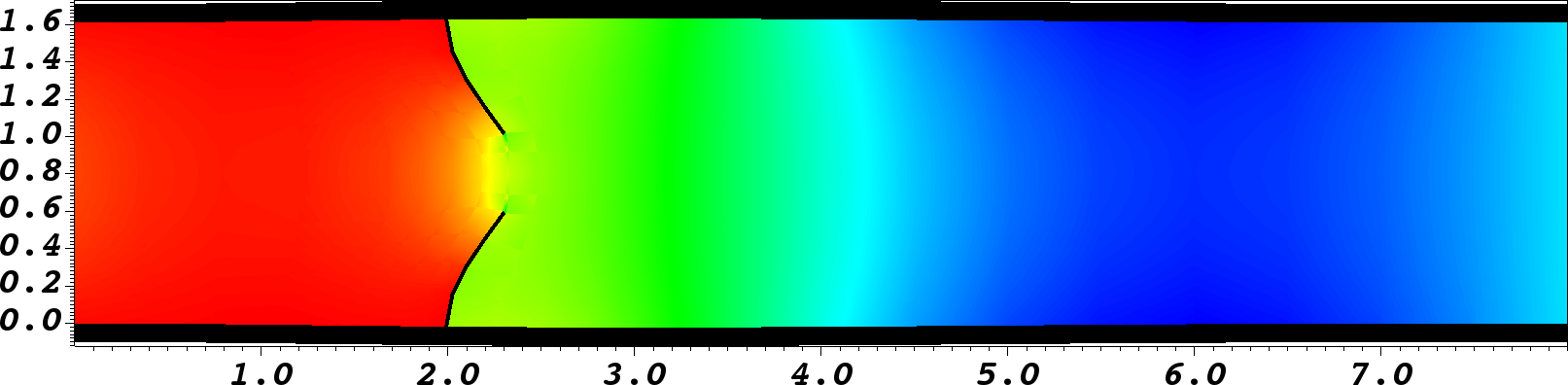}\\[1em]
  \includegraphics[width = 15cm]{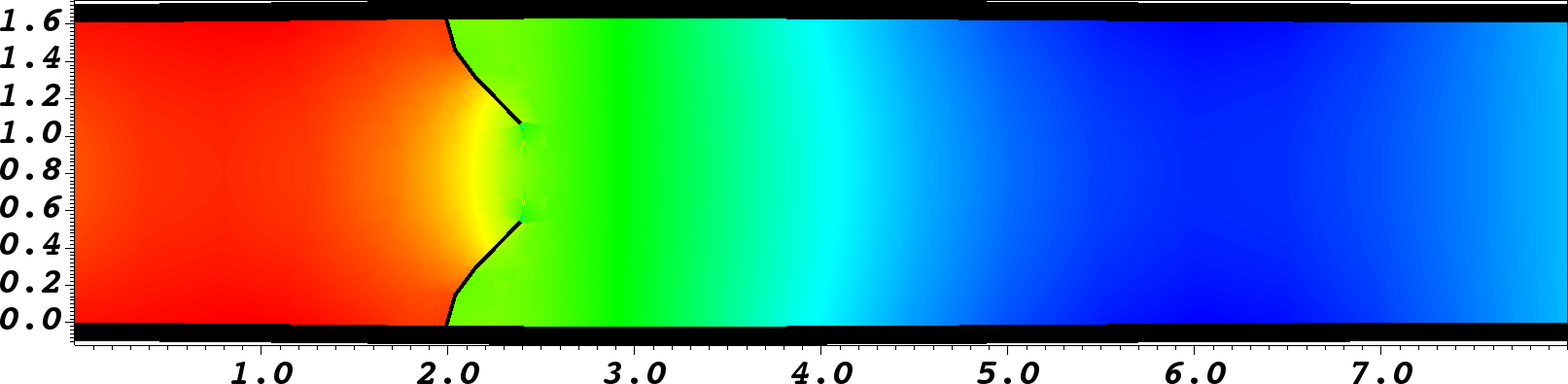}
    \caption{
      Example 3: At $T=0.579375\,\si{s}$ (time step No. $618$): $v_x(t)$ and $p(t)$ are displayed in the 
deformed configuration $\Omega(t)$. Going from top to bottom:
$v_x(t = 0.579375\,\si{s})$ in the optimization cycle $0$ (classical forward run
with $\mu = q^0 = 2\cdot 10^7$). The maximum velocity (in red) has the value 
$3.15 cm/s$. In the 2nd row, $v_x(t = 0.579375\,\si{s})$
in the eighth optimization cycle is displayed; here 
$\mu = q^8 = 5\cdot 10^6$, which means less-stiff flaps and corresponding 
higher displacements. Consequently, the maximum velocity is reduced 
and has the value $2.3\, \si{cm/s}$.
In the rows three and four the corresponding 
pressure fields are shown. The maximum pressure values 
are $3012 \,\si{\frac{g}{cm\, s^2}}$ and $2825\,\si{\frac{g}{cm\, s^2}}$, respectively.
} 
    \label{example_3_plots}
  \end{center}
 \end{figure}

The flow and pressure fields in the physical 
configuration $\Omega(t)$ are displayed in Figure~\ref{example_3_plots}.
Therein, it is visible that the solid flaps undergo large 
deformations. In the optimized configuration after $8$ cycles
the flaps even deform more. Here, a robust mesh motion 
model is indispensable. In Table~\ref{Tab:FSI-flapping}, the performance 
of the optimization procedure is shown. A reduction of 
$10^{12}$ in the cost functional is achieved. The optimal 
$q^8$ is $5\cdot 10^6$.

\newpage

\begin{table}[h!]
   \caption{Optimization results for the flapping membrane example with $\alpha = 1$ and $q_d = 5 \cdot 10^{6}$. The initial Residual in $q_0=2\cdot 10^7$ is $|\nabla
\mathcal J(q^0)| = 1.686 \cdot 10^7$}
\label{Tab:FSI-flapping}
   \centering   
   \begin{tabular}{@{}llll@{}}
\hline
Iter     & $\mathcal J(q^k)$         & $ q^k          $ &   $\frac{|\nabla
\mathcal J(q^k)|}{|\nabla
\mathcal J(q^0)|}$ \\
0	& $1.265 \cdot 10^{14}$	     & $ 2      \cdot 10^{7} $ & $ 1.0000\cdot 10^{-0}$\\     
1	& $1.9517\cdot 10^{12}$	     & $ 3.13665\cdot 10^{6} $ & $ 1.2422\cdot 10^{-1}$\\
2	& $3.0118\cdot 10^{10}$	     & $ 5.23147\cdot 10^{6} $ & $ 1.5432\cdot 10^{-2}$\\
3	& $4.6476\cdot 10^{ 8}$	     & $ 4.97125\cdot 10^{6} $ & $ 1.9170\cdot 10^{-3}$\\
4	& $7.1728\cdot 10^{ 6}$	     & $ 5.00357\cdot 10^{6} $ & $ 2.3813\cdot 10^{-4}$\\
5	& $1.1151\cdot 10^{ 5}$	     & $ 4.99956\cdot 10^{6} $ & $ 2.9582\cdot 10^{-5}$\\
6	& $2.5424\cdot 10^{ 3}$	     & $ 5.00006\cdot 10^{6} $ & $ 3.6747\cdot 10^{-6}$\\
7	& $8.6090\cdot 10^{ 2}$	     & $ 4.99999\cdot 10^{6} $ & $ 4.5649\cdot 10^{-7}$\\
8	& $8.3495\cdot 10^{ 2}$	     & $ 5      \cdot 10^{6} $ & $ 5.6707\cdot 10^{-8}$\\     
\hline
\end{tabular}
\end{table}

\newpage
\section{Conclusions}
\label{sec_conclusions}
In this work, we developed settings for FSI-based 
optimization. Therein, the FSI problem is nonlinear and nonstationary
and allows for large solid deformations. Consequently, 
when working with the ALE technique, a robust mesh motion 
model must be chosen. In this work, it is based on 
a biharmonic equation.
Based on this forward model, 
we provide the adjoint state, which is running backward-in-time.
The resulting FSI-optimization problem is solved 
with a gradient-type method.
Three numerical examples are designed to investigate the performance 
of our algorithmic techniques. In the first numerical test 
an extension of the steady-state FSI 1 benchmark is considered.
In the second and third examples, fully nonstationary tests 
are investigated. Specifically, the last 
numerical test is numerically challenging, even for the forward
problem, because the flaps 
are very thin, while undergoing large 
solid deformations. Here, we observe significant 
reductions of the cost functional and excellent
convergence properties of the optimization algorithm.





\end{document}